\documentclass[12pt]
{article}

\usepackage{amsmath,amsfonts,amssymb,amsthm}
\usepackage{latexsym}
\usepackage{makeidx}

\newtheorem{definition}{Definition}[section]
\newtheorem{theorem}[definition]{Theorem}
\newtheorem{lemma}[definition]{Lemma}
\newtheorem{corollary}[definition]{Corollary}

\newtheorem{proposition}[definition]{Proposition}

\newtheorem{question}[definition]{Question}

\newtheorem{remark}[definition]{Remark}

\newtheorem{fact}{Fact}

\newcommand\Cl{\mathrm{Cl}}

\newcommand\Zc{Zariski closure}
\newcommand\Mar{\mathfrak{M}}
\newcommand\Zar{\mathfrak{Z}}

\newcommand\ZZ{{\mathcal Z}}
\newcommand\TT{{\mathcal H}}
\newcommand\M{{\mathcal M}}

\newcommand\T{{\mathbb T}}

\newcommand\N{{\mathbb N}}

\newcommand{\topemb}{Hausdorff embedded}

\newcommand{\topembing}{Hausdorff embedding}
\newcommand{\restr}{\upharpoonright}

\title{Reflection principle characterizing groups in which unconditionally closed sets are algebraic}
\author{Dikran Dikranjan\footnote{Dipartimento di Matematica e Informatica, Unive
rsit\`{a} di Udine, Via delle Scienze  206, 33100 Udine, Italy; {\em e-mail\/}: {\tt dikranja@dimi.uniud.it}} \and Dmitri Shakhmatov\footnote{Graduate School of Science and Engineering, Division of Mathematics, Physics and Earth Sciences, Ehime University, Matsuyama 790-8577,  Japan; {\em e-mail\/}: {\tt dmitri@dpc.ehime-u.ac.jp}}}
\date{March 11, 2007}
\begin{document}
\maketitle

\begin{abstract}
We give a necessary and sufficient condition, in terms of a certain reflection principle, for every unconditionally closed subset of a group $G$ to be algebraic. As a 
corollary, we prove that this is always the case when $G$ is a 
direct product of an Abelian group with a direct product (sometimes also called a direct sum) of a family of 
countable groups. This is the widest class of groups known to date where the answer to the 
63 years old problem of Markov turns out to be positive.
We also prove that whether every unconditionally closed subset of $G$ is algebraic or not is completely determined by countable subgroups of $G$.
\end{abstract}

According to Markov \cite{M}, a subset $S$ of a group $G$ is called:
\begin{itemize}
\item[(a)] {\em elementary algebraic} if
there exist an integer $n>0$, $a_1,\ldots, a_n\in G$
and $\varepsilon_1,\ldots,\varepsilon_n\in\{-1,1\}$
such that $S=\{x\in G: x^{\varepsilon_1}a_1x^{\varepsilon_2}a_2\ldots 
a_{n-1}x^{\varepsilon_n}=a_n\}$,
\item[(b)] {\em algebraic} if $S$ is an intersection of finite unions of 
elementary algebraic subsets of $G$,
\item[(c)] {\em unconditionally closed} if $S$ is closed in
{\em every\/} Hausdorff group topology of $G$.
\end{itemize}

Since the family of all finite unions of elementary algebraic subsets 
of $G$ is closed under finite unions and contains all finite sets, it is a base 
of closed sets of some  $T_1$ topology ${\mathfrak Z}_G$ on $G$, called the {\em Zariski topology
of $G$\/}.  (This topology is also known under the name {\em verbal topology\/}, see
\cite{Bryant}.) The family of all
unconditionally closed subsets of $G$ coincides with the family of closed subsets of a $T_1$ topology 
${\mathfrak M}_G$ on $G$, namely the infimum (taken in the lattice of all topologies on $G$) of
all Hausdorff group topologies on $G$. We call ${\mathfrak M}_G$ the 
{\em  Markov topology} of $G$.  Note that $(G,{\mathfrak Z}_G)$ and $(G,{\mathfrak M}_G)$ 
are quasi-topological groups, i.e., the inversion and shifts are continuous. 

\begin{fact}
\label{Markov:stronger:than:Zariski}
 ${\mathfrak Z}_G\subseteq {\mathfrak M}_G$ for every group $G$.
\end{fact}
\begin{proof}
An elementary algebraic subset of $G$ must be closed in every Hausdorff group topology on $G$.
\end{proof}
In 1944 Markov \cite{M} asked if the equality ${\mathfrak Z}_G = {\mathfrak M}_G$ holds for every
group $G$.  He himself obtained a positive answer in case $G$ is countable: 
\begin{fact}
\label{Markov:theorem}
{\rm (Markov's theorem \cite{M})}
${\mathfrak Z}_G = {\mathfrak M}_G$ for every countable 
group $G$.
\end{fact}
Moreover, in the same manuscript \cite{M} Markov attributes to Perel'man the fact that 
${\mathfrak Z}_G = {\mathfrak M}_G$ for every Abelian group $G$. To the best of our knowledge the proof 
of this fact has never appeared in print until \cite{DS}. (We offer an alternative self-contained
proof of this result in Corollary \ref{Abelian:corollary}.)
A consistent example of a group $G$ with ${\mathfrak Z}_G\not={\mathfrak M}_G$ was announced quite recently in \cite{S1}. 

\section{Zariski and Markov embeddings} 
\label{Relations:between:embeddability}

If $H$ be a subgroup of a group $G$, then
$\Zar_G\restr_H=\{U\cap H: U\in \Zar_G\}$ denotes the subspace topology on $H$ generated by $\Zar_G$, and 
$\Mar_G\restr_H=\{U\cap H: U\in \Mar_G\}$ denotes the subspace topology on $H$ generated by 
$\Mar_G$. Note that one always have $\Zar_H\subseteq \Zar_G\restr_H$ and  $\Mar_H\subseteq \Mar_G\restr_H$.
This motivates the following definition:

\begin{definition}  {\rm 
We say that a subgroup $H$ of a group $G$ is:
\begin{itemize}
  \item[(i)] {\em Zariski embedded in $G$\/} provided that   $\Zar_H = \Zar_G\restr_H$,
i.e., the subspace topology induced on $H$ by the Zariski topology of $G$ coincides with the Zariski topology of $H$, 
  \item[(ii)] {\em Markov embedded in $G$\/} provided that  $\Mar_H = \Mar_G\restr_H$,
i.e.,  the Markov topology of $H$  coincides with the subspace topology induced on $H$ by the Markov topology of $G$.
\end{itemize}
}
\end{definition}

We shall see in the sequel that every subgroup $H$ of an Abelian group 
$G$ is both Zariski embedded and Markov embedded in $G$. 
For every infinite Abelian group $H$ there exists a (necessarily non-Abelian)
group $G$ containing $H$ as a subgroup such that $H$ is not Markov embedded in $G$ (see Remark 
\ref{absolutely-Markov:implies:absolutely-Hausdorff?}(ii)).
An example of a subgroup $H$ of a (necessarily non-Abelian) group 
$G$ that is neither Zariski embedded nor Markov embedded in $G$ 
can be found in 
Remark \ref{Example:nonZar/Mar_embedded}. 

Distinguishing Zariski and Markov embeddings is surprisingly difficult. 
Indeed, our next lemma indicates that the difference, if any, is closely related to Markov's problem.

\begin{lemma}
\label{Markov:vs:Zariski:embedding} Let $H$ be a subgroup of a group $G$.
\begin{itemize}
\item[(a)] if $\Zar_H=\Mar_H$ and $H$ is Markov embedded in $G$, then $H$ is also Zariski embedded in $G$. 
\item[(b)] if $\Zar_G=\Mar_G$ and $H$ is  Zariski embedded in $G$, then $H$ is also Markov embedded in $G$.
\end{itemize}
\end{lemma}
\begin{proof}
(a) From  Fact \ref{Markov:stronger:than:Zariski}, we have $\Zar_G\restr_H\subseteq \Mar_G\restr_H$.
Since $H$ is Markov embedded in $G$, we have $\Zar_H\subseteq \Zar_G\restr_H\subseteq \Mar_G\restr_H=\Mar_H$.
Since $\Zar_H =\Mar_H$, we get $\Zar_H= \Zar_G\restr_H$. 
This means that $H$ is Zariski embedded in $G$.

(b) From the assumptions of (b) we get  $\Zar_H= \Zar_G\restr_H= \Mar_G\restr_H$. Since $ \Zar_H\subseteq \Mar_H$,
this proves $\Mar_G\restr_H\subseteq \Mar_H$. Since the converse inclusion $\Mar_H\subseteq \Mar_G\restr_H$
always holds, we obtain $\Mar_H = \Mar_G\restr_H$. This means that $H$ is Markov embedded in $G$.
\end{proof}

\begin{remark}\label{Remark:ZvsM}
{\rm 
A careful analysis of the above proof  reveals that $\Zar_H=\Mar_H$ holds under the assumption of item (b) of Lemma \ref{Markov:vs:Zariski:embedding}.
}
\end{remark}

\begin{corollary}
\label{Markov:embedding:coincides:with:Zariski:embedding:for;countable:groups}
Let $H$ be a countable subgroup of a group $G$.
\begin{itemize}
\item[(a)]  If $H$ is Markov embedded in $G$, then $H$ is also Zariski embedded in $G$. 
\item[(b)]  If, in addition, also $G$ is countable, then $H$ is Markov embedded in $G$ if and only if $H$ is Zariski embedded in $G$. 
\end{itemize}
\end{corollary}

\begin{proof}
Immediately follows from Lemma \ref{Markov:vs:Zariski:embedding} and Fact \ref{Markov:theorem}.
\end{proof}

\begin{corollary}\label{LAST:corollary}
\label{Markov:embedding:coincides:with:Zariski:embedding:for_countable:groups}
Let $H$ be a Markov embedded subgroup of a group $G$ with 
$\Mar_G=\Zar_G$. Then $H$ is Zariski  embedded in $G$ if and only if $\Mar_H=\Zar_H$.\end{corollary}

\begin{proof} Apply Remark \ref{Remark:ZvsM} and Lemma \ref{Markov:vs:Zariski:embedding}.
\end{proof}

\section{\topembing s}

\begin{definition}
{\rm \cite{DS1}
A subgroup $H$ of a group $G$ is called:

\begin{itemize}
  \item[(i)]  {\em \topemb\  in $G$} provided that every Hausdorff group topology $\mathcal{T}$ on $H$ is a restriction of some Hausdorff group topology $\mathcal{T}^*$ on $G$ (and in this case we say that $\mathcal{T}^*$ {\em extends $\mathcal{T}$\/}),
\item[(ii)] {\em super-normal (in $G$)} provided that for every $x \in G$ there exists $y \in H$ such that $x^{-1}hx=y^{-1}hy$ for all $h \in H$. 
\end{itemize}
}
\end{definition}

Obviously, super-normal subgroups are normal. 

\begin{lemma}\label{(b)} A normal subgroup $H$ of a group $G$ is super-normal if and only if 
$G= c_G(H) H$, where $c_G(H)$ is the centralizer of $H$ in $G$. 
\end{lemma}
\begin{proof}
Assume $H$ is super-normal and pick an element $x\in G$. Then there exists $y\in H$ such that $y^{-1}hy=x^{-1}hx$ 
for all $h \in H$. Then $xy^{-1}\in c_G(H)$, so $x\in c_G(H) H$. 

If $G=c_G(H) H$, then for $x\in G$ there exists $y\in H$ such that 
$x\in c_G(H) y$. This means that $xy^{-1}\in c_G(H)$, and hence
$xy^{-1}h=h xy^{-1}$ for all $h \in H$. This yields $y^{-1}hy=x^{-1}hx$  for every $h \in H$.
\end{proof}

The lemma gives the following immediate corollary. 
 
\begin{corollary}\label{Cor(b)} Every direct summand, as well as every 
central subgroup, is super-normal. In particular, every subgroup of an Abelian group is super-normal.
\end{corollary}

The next theorem characterizing \topemb\ normal subgroups is taken from {\rm \cite{DS1}}. We give its proof here for the reader's  convenience. 

\begin{theorem} 
\label{characterization:of:normal:topologically.embedded.subgroups}
Let $N$ be a normal subgroup of the group $G$. Then 
$N$ is \topemb\  in $G$ iff the automorphisms of $N$ induced by conjugation by elements of  $G$ are continuous  for any group topology on $N$.
\end{theorem}

\begin{proof}
The necessity is obvious since the conjugations are continuous in any  topological group. Assume now that all automorphisms of $N$ induced by 
the conjugation by elements of $G$ are $\mathcal{T}$-continuous for any  Hausdorff group topology $\mathcal{T}$ on $N$. 
Fix a Hausdorff group topology $\mathcal{T}$ on $N$. Take the filter of all neighbourhoods of 1 in $(N,\mathcal{T})$ as a 
base of neighbourhoods of 1 in a group topology $\sigma$ of $G$. This works 
since the only axiom to check is to find, for every $x\in G$ and every $\sigma$-neighbourhood $U$ of 1, a $\sigma$-neighbourhood $V$ of 1 such that 
$V^x:=x^{-1}Vx\subseteq U$. Since we can choose $U,V$ contained in 
$N$, this immediately follows from our assumption of $\mathcal{T}$-continuity of  the restrictions to $N$ of the conjugations in $G$.
\end{proof}

\begin{corollary}\label{cyclic:subgroup}
Every normal  cyclic subgroup is Hausdorff embedded.  
\end{corollary}

\begin{proof}
Assume $H$ is a normal cyclic subgroup of a group $G$. Then every automorphism of $H$ is continuous in any group topology of $H$.
Therefore, Theorem \ref{characterization:of:normal:topologically.embedded.subgroups} applies. 
\end{proof} 

\begin{corollary}
\label{super-notmal:implies:topemb}
If a subgroup $H$ of a group $G$ is super-normal in $G$, then $H$ is  \topemb\ in $G$.
\end{corollary}

\begin{proof}
As $H$ is super-normal, each conjugation by an element  of $G$ coincides with the conjugation by some element of $H$, so each such conjugation is 
continuous  in any group topology on $H$. Now Theorem \ref{characterization:of:normal:topologically.embedded.subgroups} applies.  
\end{proof}

The implication in the above corollary is not reversible:
a normal \topemb\ subgroup $H$ of a group $G$ need not be super-normal in $G$ \cite{DS1}. 

The following lemma is obvious.

\begin{lemma}
\label{topological:embedding:implies:markov:embedding}\label{super-normal:theorem}
Let $H$ be a subgroup of a group $G$. If $H$ is \topemb\ in $G$, then $H$ is also Markov embedded in $G$.
\end{lemma}

In \cite{DS1}, a normal subgroup of a countable
group $G$ is constructed such that $H$ is Zariski embedded in $G$ but not \topemb\ in $G$. By Corollary \ref{Markov:embedding:coincides:with:Zariski:embedding:for;countable:groups}(b), $H$ is also Markov embedded in $G$.
This shows that the implication of Lemma \ref{topological:embedding:implies:markov:embedding} is not reversible, even for a normal subgroup $H$.

\begin{remark}\label{former:item:(e)}
{\rm If $h:G\to G_1$ is a group isomorphism and the  subgroup $H$ of $G$ is \topemb\ (Markov embedded, Zariski embedded) in $G$, then 
   the subgroup $f(H)$ of $G_1$
is \topemb\ (respectively, Markov embedded, Zariski embedded)  in $G_1$.  
}
\end{remark}

Our next result uncovers a curious fact: If a countable subgroup $H$ of a group 
$G$ fails to be \topemb, then this failure can always be witnessed by some 
{\it metric\/} group topology on $H$.

\begin{theorem}\label{metric:groups}
Let $H$ be a countable subgroup of a group $G$. If every metric group topology on $H$ can be extended to a (not necessarily metric) group topology on $G$, then $H$ is \topemb\ in $G$.  
\end{theorem}
\begin{proof}
Let $\mathcal{T}$ be a Hausdorff group topology on $H$. Then $\mathcal{T}$ has a countable network, and the main result of \cite{Shakh} implies that $\mathcal{T}$ is the supremum of some family $\{\mathcal{T}_i:i\in I\}$ of group topologies on $H$  with a countable base. Then each $\mathcal{T}_i$ is metric, and so by the assumption of our lemma, there exists a Hausdorff group topology $\mathcal{T}^*_i$ on $G$ extending $\mathcal{T}_i$. Now the supremum of the family 
 $\{\mathcal{T}^*_i:i\in I\}$ is the Hausdorff group topology on $G$ that obviously induces
$\mathcal{T}$ on $H$. Hence, $H$ is \topemb.
\end{proof}

Let $\M$ ($\ZZ,\TT$) denote the class of group embeddings $H\hookrightarrow G$
such that $H$ is Markov (resp., Zariski, Hausdorff) embedded in $G$. Then 
one can easily verify that $\M$, $\ZZ$ and $\TT$  are stable under composition and left cancellation. More precisely:

\begin{lemma}
\label{general:properties:of:MZ:embeddings}
 If  $H_1 \leq H_2\leq G$ are groups and $i_1:H_1\hookrightarrow H_2$, $i_2:H_2\hookrightarrow G$ 
and $i_2\circ i_1: H_1\hookrightarrow G$ are the respective inclusions, then: 
\begin{itemize}
  \item[(a)] $i_1, i_2 \in \M$ ($i_1, i_2 \in \ZZ, \TT$) implies $i_2\circ  i_1\in \M$ (respectively, $i_2\circ  i_1\in \ZZ$, $i_2\circ  i_1\in \TT$);
  \item[(b)] if $i_2\circ i_1 \in \M$ ($i_2\circ i_1 \in \ZZ, \TT$) , then also $i_1\in   \M$ (respectively, $i_1\in   \ZZ$, $i_1\in   \TT$).
  \end{itemize}
\end{lemma} 

\begin{lemma}\label{(*)} 
Let $H$ be a subgroup of the direct product $G=G_1\times G_2$. If the subgroup $G_1\cap H$ of $G_1$ 
is Markov embedded in $G_1$, then $G_1\cap H$ is Markov embedded in $H$ as well.
\end{lemma}
\begin{proof}
 Indeed,  as a direct summand of $G$, $G_1$ is super-normal in $G$. By Corollary
\ref{super-notmal:implies:topemb} and Lemma
\ref{topological:embedding:implies:markov:embedding}, 
$G_1$ is Markov embedded in $G$, so Lemma \ref{general:properties:of:MZ:embeddings}(a) allows us to conclude that $G_1\cap H$ is Markov embedded in $G$. Applying now Lemma \ref{general:properties:of:MZ:embeddings}(b) we conclude that 
$G_1\cap H$ is Markov embedded in $H$.
\end{proof}

By Lemma \ref{general:properties:of:MZ:embeddings}, the classes $\M$ and $\TT$ (of Markov embeddings and 
\topembing s) are closed under composition and left cancellation. Now we  
are going to show that these classes are not closed under pullback. More precisely, if $G_1 \hookrightarrow  G$ is a 
\topembing\  and $H$ is a subgroup of $G$, then we shall see 
that the induced embedding $G_1\cap H \hookrightarrow H$  need not be even a Markov embedding. We  take $G$ of the special form $G=G_1 \times G_2$, so that $G_1 \hookrightarrow  G$, being the inclusion of a direct summand, is certainly a 
\topembing, hence a Markov embedding (Lemma \ref{topological:embedding:implies:markov:embedding}). Then for an appropriate subgroup $H$ of $G=G_1 \times G_2$ 
 we show that $G_1\cap H\hookrightarrow H$ is not a even a Markov embedding. By Lemma \ref{(*)}, this will show that also $G_1\cap H\hookrightarrow G_1$ fails to be a Markov embedding.

\begin{lemma}
\label{lemma:of:product:of:two}
Let $N$ be a countable Abelian group that admits a decomposition $N=N_1 
\times N_2$ into a direct product of two infinite groups $N_1$ and $N_2$. 
 Then there exist a countable group $G'$, a subgroup $H$ of the direct product $G=G_1\times G_2$, where $G_1=G_2=G'$, and a metric group topology
$\mathcal{T^*}$ on $G^*=H\cap G_1$ having the following properties:
\begin{enumerate}
  \item[(i)] $G^*$ is isomorphic to $N$ and $[H:G^*]=2$,
  \item[(ii)]  $G^*$ is neither Markov nor Zariski embedded in $H$,
  \item[(iii)] $\mathcal{T^*}$ cannot be extended to any Hausdorff group topology on $H$.
\end{enumerate}
\end{lemma}
 
 \begin{proof} By \cite[Lemma 3.9]{DS1}, there exists an involution $f$ such that  $N$ is not Zariski embedded in the 
 countable semidirect product $G'=N\rtimes\langle f\rangle$. 
(Here $\langle f\rangle$ denotes the two-element cyclic group generated by the involution $f$.)
Since $G'$ is countable, $N$ is not Markov embedded in $G'$ by Corollary
\ref{Markov:embedding:coincides:with:Zariski:embedding:for;countable:groups}.
Define $G_1=G_2=G'$, and let $H$ be the subgroup of $G=G_1\times G_2$ generated by the element $(f,f)\in G$ 
and the subgroup $N\times \{1_{G_2}\}$ of $G$.
Note that the projection $p_1:G=G_1\times G_2\to G_1$ onto the first coordinate sends 
$H$ isomorphically onto $G_1=G'$ and $p_1(G^*)=N$. This proves (i). 

Since $N$ is neither Markov nor Zariski embedded in  $G'$ and $p_1$ sends $H$ isomorphically onto $G'$ with $p_1(G^*)=N$, it follows from Remark \ref{former:item:(e)} that $G^*$ is neither Markov nor Zariski embedded in $H$. This proves (ii). 

To prove (iii), note that $G^*$ is not \topemb\ in $H$ by item (ii) and Lemma \ref{topological:embedding:implies:markov:embedding}. Since $G^*$ is countable, 
 by Theorem \ref{metric:groups}
there must exist a metric group topology $\mathcal{T}^*$ on $G^*$ that cannot be extended to any Hausdorff group topology on $H$.
\end{proof}

\section{Reflection principle for the Zariski closure operator}\label{Refl}

If $X$ is a set, then $[X]^{<\omega}$ and $[X]^{\le\omega}$ denote the set of all finite subsets of $X$ and all (at most) countable subsets of $X$, respectively. 
$\mathbb{N}$ denotes the set of all natural numbers. We need some machinery from set theory useful for carrying out closing off arguments.

\begin{definition} {\rm 
Let  $\mathcal{C}\subseteq [X]^{\le\omega}$.

\begin{itemize}
   \item[(i)] $\mathcal{C}$  is {\em closed in $[X]^{\le\omega}$\/}  if, whenever $\{C_n:n\in\mathbb{N}\}\subseteq \mathcal{C}$ 
and $C_0\subseteq C_1\subseteq \dots\subseteq C_n\subseteq C_{n+1}\subseteq\dots$, then $\bigcup\{C_n:n\in\mathbb{N}\}\in \mathcal{C}$,

\item[(ii)] $\mathcal{C}$  is {\em unbounded in $[X]^{\le\omega}$\/} provided that for every $Y\in [X]^{\le\omega}$ there exists $C\in\mathcal{C}$ with 
$Y\subseteq C$,

\item[(iii)] $\mathcal{C}$  is a {\em club in $[X]^{\le\omega}$\/} (a common abbreviation for ``closed and unbounded'') if $\mathcal{C}$ is both closed and unbounded in $[X]^{\le\omega}$. 
\end{itemize}
}
\end{definition}

For a group $G$ we define $\mathcal{S}(G)=\{H\in[G]^{\le\omega}: H$ is a subgroup of $G\}$. This is a typical example of a club:

\begin{lemma}
\label{subgroup:lemma}
If $G$ is a group, then $\mathcal{S}(G)$ is a club in $[G]^{\le\omega}$.
\end{lemma}

As witnessed by (the proof of) the previous lemma, clubs appear naturally in various closing off arguments, and a general scheme 
that greatly simplifies carrying out such arguments is given below.

\begin{definition}
Given a set $X$ and a function $\varphi:[X]^{<\omega}\to[X]^{<\omega}$, we say that
a subset $Y$ of $X$ is {\em $\varphi$-invariant\/} provided that  $\varphi([Y]^{<\omega})\subseteq [Y]^{<\omega}$.
\end{definition}

\begin{lemma}
\label{closure:lemma}
Given a set $X$ and a function $\varphi:[X]^{<\omega}\to[X]^{<\omega}$, the family
$\mathcal{I}(\varphi)=\{Y\in[X]^{\le\omega}: Y$ is $\varphi$-invariant$\}$ is a club in $[X]^{\le\omega}$.
\end{lemma}

\begin{proof}
One can easily check that $\mathcal{I}(\varphi)$ is closed in $[X]^{\le\omega}$. Let  us show that $\mathcal{I}(\varphi)$ is also unbounded in $[X]^{\le\omega}$.
Fix arbitrarily $Y\in [X]^{\le\omega}$. By induction on $n\in\mathbb{N}$
define a sequence $\{Y_n:n\in\mathbb{N}\}\subseteq [X]^{\le\omega}$ by $Y_0=Y$
and $Y_{n+1}=\bigcup\{\varphi(Z):Z\in[Y_n]^{<\omega}\}\cup Y_n$. Finally, note that $E=\bigcup\{Y_n:n\in\mathbb{N}\}$
is $\varphi$-invariant and $Y\subseteq E$.
\end{proof}

The following well-known lemma reveals one of the main reasons why clubs are so useful. We briefly outline the proof for the reader's convenience.

\begin{lemma}
\label{intersection:of:clubs}
If $\{\mathcal{C}_n:n\in\mathbb{N}\}$ is a sequence of clubs in $[X]^{\le\omega}$, then 
$\mathcal{C}=\bigcap\{\mathcal{C}_n:n\in\mathbb{N}\}$ is also a club in $[X]^{\le\omega}$.
\end{lemma}

\begin{proof}
Clearly, $\mathcal{C}$ is closed in $[X]^{\le\omega}$. Let us show that $\mathcal{C}$ is also unbounded in $[X]^{\le\omega}$.
Fix arbitrarily $Y\in [X]^{\le\omega}$. Since each $\mathcal{C}_n$ is unbounded 
in $[X]^{\le\omega}$, there exists a function $f_n:[X]^{\le\omega}\to \mathcal{C}_n$ such that $Z\subseteq f_n(Z)$ for all 
$Z\in [X]^{\le\omega}$. Fix an enumeration $\mathbb{N}\times\mathbb{N}=\{(k_n,m_n):n\in\mathbb{N}\}$ of $\mathbb{N}\times\mathbb{N}$.  By induction on $n\in\mathbb{N}$
define a sequence $\{Y_n:n\in\mathbb{N}\}\subseteq [X]^{\le\omega}$ by $Y_0=Y$
and $Y_{n+1}=Y_n\cup f_{k_n}(Y_n)$.
Then $C=\bigcup\{Y_n:n\in\mathbb{N}\}\in\mathcal{C}$ and $Y\subseteq C$.
\end{proof}

Let $G$ be a group. Given $n\in\mathbb{N}$, $a\in  G^{n+1}$ and $\varepsilon\in \{-1,1\}^{n+1}$ we define
$$
E_n({a,\varepsilon}; G)=\{x\in G: x^{\varepsilon(0)}a(0)x^{\varepsilon(1)}a(1)\dots
a(n-1)x^{\varepsilon(n)}=a(n)
\}
$$
and $S_n(a,\varepsilon)=\{a(0), a(1),\dots,a(n)\}$. Define 
$$\mathcal{F}_G=\bigcup_{n\in\mathbb{N}}\{n\}\times G^{n+1}\times \{-1,1\}^{n+1},$$
and for $F\in[\mathcal{F}_G]^{<\omega}\setminus\{\emptyset\}$ 
 let $S(F)=\bigcup \{S_n(a,\varepsilon): (n,a,\varepsilon)\in F\}$.
Define also $S(\emptyset)=\emptyset$.

For $F\in [\mathcal{F}_G]^{<\omega}$ define 
$U_G(F)=G\setminus \bigcup_{(n,a,\varepsilon)\in F} E_n(a,\varepsilon;G)$. Clearly, the family 
$\{U_G(F):F\in [\mathcal{F}_G]^{<\omega}\}$ forms a base of the Zariski topology $\Zar_G$ on $G$, and 
the closure 
$$
\Cl_{\Zar_G}A=\{z\in G\ :\  \forall\ F\in [\mathcal{F}_G]^{<\omega}\ 
(z\in U_G(F)\to A\cap U_G(F)\not=\emptyset)\}
$$
of a set $A\subseteq G$ in this topology
 is called the {\em \Zc\ of $A$ in $G$\/}.

The main result of this section is the following general reflection principle for the Zariski closure.
\begin{theorem}
\label{reflection:lemma}
Let $G$ be a group and $A$ a subset of $G$. Then the family 
$$
\mathcal{Z}_A=\{H\in \mathcal{S}(G):\Cl_{\Zar_H}(H\cap A)= H\cap \Cl_{\Zar_G}A\}
$$
contains a club in $[G]^{\le\omega}$.
\end{theorem}

\begin{proof}
For every $F\in [\mathcal{F}_G]^{<\omega}$, if $A\cap U_G(F)\not=\emptyset$,  pick some
\begin{equation}
\label{equation.2}
x_{F}\in A\cap U_G(F),
\end{equation}
and define $x_{F}=e$ otherwise. (Here $e$ denotes the identity element of $G$.) 

For $z\in  \Cl_{\Zar_G}A$ define $F_{z}=\emptyset$.
For  $z\in G\setminus \Cl_{\Zar_G}A$, choose $F_{z}\in [\mathcal{F}_G]^{<\omega}$ satisfying
\begin{equation}
\label{equation.1}
A\cap U_G(F_{z})=\emptyset \mbox{ and }
z\in U_G(F_{z}).
\end{equation}

Define functions
$\varphi_{k}: [G]^{<\omega}\to [G]^{<\omega}$ 
(for every $k\in\mathbb{N}$) and $\psi:[G]^{<\omega}\to [G]^{<\omega}$ by 
$$
\varphi_{k}(X)=\left\{x_{F}\ \left|\ \  F\subseteq \bigcup_{n\le k} \{n\}\times  X^{n+1}\times \{-1,1\}^{n+1}\right\}\right.
$$
and $\psi(X)=\bigcup \{S(F_{z}) :{z\in X}\}$ for $X\in [G]^{<\omega}$. 

According to Lemmas \ref{subgroup:lemma}, \ref{closure:lemma}  and
\ref{intersection:of:clubs}, the family $\mathcal{H}=\mathcal{S}(G)\cap\mathcal{I}(\psi)\cap
\bigcap_{k\in\mathbb{N}}\mathcal{I}(\varphi_{k})$ is a club in $[G]^{<\omega}$. 
It remains only to show that  $\mathcal{H}\subseteq \mathcal{Z}_A$.

Fix $H\in \mathcal{H}$. We have to check that $\Cl_{\Zar_H}(H\cap A)= H\cap \Cl_{\Zar_G}A$. To start with, note that 
\begin{equation}
\label{equation.5}
U_H(F)=H\cap U_G(F)\mbox{ for every } F\in[\mathcal{F}_H]^{<\omega}.
\end{equation}

First, let us show that 
$\Cl_{\Zar_H} (H\cap A)\subseteq H\cap \Cl_{\Zar_G}A$.
Pick arbitrarily $z\in H\setminus  \Cl_{\Zar_G}A\subseteq G\setminus  \Cl_{\Zar_G}A$. 
By our choice of $F_{z}$, (\ref{equation.1}) holds.
Note that $S(F_{z})\subseteq \psi(\{z\})\in[H]^{<\omega}$ because $z\in H$ and $H\in\mathcal{I}(\psi)$.
Therefore, $S(F_{z})\in [H]^{<\omega}$, which in turn yields
$F_{z}\in [\mathcal{F}_H]^{<\omega}$. From (\ref{equation.5}), (\ref{equation.1}) and $z\in H$ we get
$z\in H\cap U_G(F_{z})=U_H(F_{z})$ and 
$
(H\cap A)\cap U_H(F_{z}) =
(H\cap A)\cap H\cap U_G(F_{z})\subseteq A\cap U_G(F_{z})=\emptyset.
$
This yields $z\not\in\Cl_{\Zar_H}(H\cap A)$.

Second, let us prove the inverse inclusion $H\cap \Cl_{\Zar_G}A\subseteq \Cl_{\Zar_H} (H\cap A)$.
Pick arbitrarily $z\in H\cap \Cl_{\Zar_G}A$.  Assume that $F\in[\mathcal{F}_H]^{<\omega}$ and  $z\in U_H(F)$.
We are going to show that $(H\cap A)\cap U_H(F)\not=\emptyset$.
From $z\in H$ and (\ref{equation.5}) it now follows that  $z\in U_G(F)$.
From
$F\in[\mathcal{F}_H]^{<\omega}\subseteq [\mathcal{F}_G]^{<\omega}$
and $z\in\Cl_{\Zar_G}A$ we must also have 
$A\cap U_G(F)\not=\emptyset$, and thus (\ref{equation.2})  holds by our choice of $x_{F}$.

Let $k=\max\{n\in\mathbb{N}:(n,a,\varepsilon)\in F\}$. Then  $x_{F}\in \varphi_{k}(S(F))$.
From $F\in[\mathcal{F}_H]^{<\omega}$, it follows that  $S(F)\in [H]^{<\omega}$.
Since $H\in\mathcal{I}(\varphi_{k})$, $H$ is $\varphi_{k}$-invariant,
and thus $\varphi_{k}(S(F))\in [H]^{<\omega}$. We conclude that $x_{F}\in H$.
Combining this with (\ref{equation.2}) and (\ref{equation.5}), we get
$$
x_{F}\in H\cap (A\cap U_G(F))=(H\cap A)\cap (H\cap U_G(F))
=(H\cap A)\cap U_H(F)\not=\emptyset.
$$
\end{proof}

From Theorem \ref{reflection:lemma} and Lemma \ref{intersection:of:clubs}, we obtain the following
\begin{corollary}
\label{reflection:corollary}
Let $G$ be a group and $\mathcal{A}$ a countable family of subsets of $G$.
Then the family
$\{H\in \mathcal{S}(G):\Cl_{\Zar_H}(H\cap A)= H\cap \Cl_{\Zar_G}A\mbox{ for all }A\in\mathcal{A}\}$
contains a club in $[G]^{\le\omega}$.
\end{corollary}

\begin{corollary}
Let $G$ be a group and $\mathcal{A}$ a countable family of $\Zar_G$-closed subsets of $G$. 
Then the family 
$\{H\in \mathcal{S}(G):H\cap A\mbox{ is $\Zar_H$-closed for every }A\in\mathcal{A}\}$
contains a club in $[G]^{\le\omega}$.
\end{corollary}

\begin{corollary}
Assume that $G$ is a group, $\mathcal{A}$ a countable family of $\Zar_G$-closed subsets of $G$ and $X$ is a countable subset of $G$. Then there exists a countable subgroup $H$ of $G$ containing $X$ such that
$H\cap A$ is $\Zar_H$-closed for each $A\in\mathcal{A}$.
 \end{corollary}

\begin{remark}
{\rm 
For a reader familiar with the notion of elementary submodels we note in passing that an alternative proof of Theorem \ref{reflection:lemma} could be furnished using model-theoretic methods. Indeed, the family $\mathcal{C}$ consisting 
of all intersections $M\cap G$, where $M$ is a countable elementary submodel of (sufficiently large fragment) of the universe containing $(G,\cdot,{}^{-1})$ and $A$, forms a club
in $[G]^{\le\omega}$ satisfying $\mathcal{C}\subseteq \mathcal{Z}_A$.
}
\end{remark}

\section{Characterization of groups for which Markov and Zariski topologies coincide}\label{applications:section}

It turns out that the version of reflection for $\Mar_G$ similar to the one for $\Zar_G$ obtained in 
Theorem \ref{reflection:lemma} {\em characterizes\/} groups $G$ for which Markov and Zariski topologies coincide. 

\begin{theorem}
\label{solving:Markov}
For a group $G$ the following conditions are equivalent:

\begin{itemize}
   \item[(i)]  $\Mar_G=\Zar_G$;

   \item[(ii)] For every set $A\subseteq G$, the family
$$\mathcal{M}_A=\{H\in \mathcal{S}(G):\Cl_{\Mar_H}(H\cap A)= H\cap \Cl_{\Mar_G}A\}$$
contains a club in $[G]^{\le\omega}$;

   \item[(iii)] For every $\Mar_G$-closed set $A\subseteq G$, the family
$\mathcal{E}_A=\{H\in\mathcal{S}(G):H\cap A\mbox{ is $\Mar_H$-closed}\}$ 
contains a club in $[G]^{\le\omega}$.
\end{itemize}
\end{theorem}
\begin{proof}
(i)$\to$(ii).
Let $A$ be a subset of $G$.
Applying Theorem \ref{reflection:lemma}, we conclude that the family $\mathcal{Z}_A$ (as given by Theorem \ref{reflection:lemma}) contains some club in $[G]^{\le\omega}$. Therefore, it suffices to show that
$\mathcal{Z}_A\subseteq\mathcal{M}_A$. Let $H\in \mathcal{Z}_A$. From (i) and the definition of $\mathcal{Z}_A$, we get
$\Cl_{\Zar_H}(H\cap A)= H\cap \Cl_{\Zar_G}A=H\cap \Cl_{\Mar_G} A$.
Since $H$ is countable, from  Fact \ref{Markov:theorem}
it follows that $\Cl_{\Mar_H}(H\cap A)=\Cl_{\Zar_H}(H\cap A)=H\cap \Cl_{\Mar_G} A$,
which yields $H\in \mathcal{M}_A$.

(ii)$\to$(iii) is trivial.

(iii)$\to$(i). We have to show that every $\Mar_G$-closed set is $\Zar_G$-closed.
Suppose that some $\Mar_G$-closed set $A$ is not $\Zar_G$-closed. Then there exists $g\in \Cl_{\Zar_G}A\setminus A$. 
Let $\mathcal{Z}_A$ be the family from the conclusion of Theorem \ref{reflection:lemma}. 
By (ii), Theorem  \ref{reflection:lemma} and  Lemma \ref{intersection:of:clubs}, $\mathcal{E}_A\cap \mathcal{Z}_A$
contains some club $\mathcal{C}$. Since $\mathcal{C}$ is unbounded, there
exists $H\in\mathcal{C}$ with $g\in H$. From $H\in \mathcal{E}_A$, it follows that $H\cap A$ is $\Mar_H$-closed.
Since $H$ is a countable group, $H\cap A$ must also be $\Zar_H$-closed by Fact \ref{Markov:theorem}.
Since $H\in \mathcal{Z}_A$, we have $g\in H\cap \Cl_{\Zar_G}A=\Cl_{\Zar_H}(H\cap A)= H\cap A$, 
in contradiction with $g\not\in A$.
\end{proof}

\begin{corollary}
The equality $\Mar_G=\Zar_G$ is completely determined by countable subgroups of $G$.
\end{corollary}

As an application of Theorem \ref{solving:Markov}, we get a new class of groups $G$ for which ${\mathfrak Z}_G={\mathfrak M}_G$:

\begin{corollary}
\label{approximation:theorem}
Let $G$ be a group such that the family $\mathcal{N}_G=\{N\in \mathcal{S}(G): N$ is Markov embedded in $G\}$ 
contains some club in $[G]^{\le\omega}$. Then Markov and Zariski topologies on $G$ coincide.
\end{corollary}

\begin{proof}
Indeed, given $\Mar_G$-closed set $A\subseteq G$, we have $\mathcal{N}_G\subseteq \mathcal{M}_A$, and the conclusion follows from the implication
(iii)$\to$(i) of Theorem \ref{solving:Markov}.
\end{proof}

\begin{corollary}
{\rm \cite{DS}}
\label{Abelian:corollary}
Markov and Zariski topologies coincide for Abelian groups.
\end{corollary}

\begin{proof} Let $H$ be a subgroup of an Abelian group $G$. By Corollary \ref{Cor(b)}, $H$ is super-normal
in $G$. Hence, by Corollary \ref{super-notmal:implies:topemb} and Lemma \ref{super-normal:theorem}, $H$ is Markov embedded in $G$. 
We have proved that $\mathcal{S}(G)\subseteq \mathcal{N}_G$. Now the conclusion of our corollary follows from Lemma
\ref{subgroup:lemma} and Corollary  \ref{approximation:theorem}.
\end{proof}

Markov \cite{M} has attributed (the equivalent form of) Corollary  \ref{Abelian:corollary} to  Perel'man. To the best of 
our knowledge the proof has never appeared in print until \cite{DS}.  (In the particular case when $G$ is almost torsion-free\footnote{An Abelian group $G$ is {\em 
 almost torsion-free} if $G[n]=\{g\in G: ng=0\}$ is finite for every $n>1$.}  the equality ${\mathfrak 
Z}_G={\mathfrak M}_G$  was earlier proved in \cite{TY}.) In fact, \cite{DS} also offers a much stronger version of this result.

Our next result is a counterpart of Corollary 
\ref{Markov:embedding:coincides:with:Zariski:embedding:for;countable:groups}(a).
\begin{corollary}\label{Mar:Embedded>Zariski:Embedded:corollary}
Let $H$ be  an Abelian subgroup of a group $G$. If $H$ is Markov embedded in $G$, then $H$ is also Zariski embedded in $G$. 
\end{corollary}
\begin{proof} $\Zar_H=\Mar_H$ by Corollary \ref{Abelian:corollary}.
Now Lemma \ref{Markov:vs:Zariski:embedding} (a) applies. 
\end{proof}

Let $\{G_i:i\in I\}$ be a family of groups. We denote by $\bigoplus_{i\in I} G_i$
the set of all functions $g: I\to \bigcup_{i\in I} G_i$ such that 
$g(i)\in G_i$ for all $i\in I$ and the set $\{i\in I: g(i)\not= 1_i\}$ is finite.
(Here $1_i$ denotes the identity element of $G_i$.) 
For $g,h\in \bigoplus_{i\in I} G_i$
define functions $gh: I\to \bigcup_{i\in I} G_i$ 
and $g^{-1}:I\to \bigcup_{i\in I} G_i$
by $gh(i)=g(i)h(i)$ and $g^{-1}(i)=(g(i))^{-1}$ for all $i\in I$.
It is easy to check that with these two operations $\bigoplus_{i\in I} G_i$
becomes a group which 
we will call the 
{\em direct sum\/} 
of the family 
$\{G_i:i\in I\}$. While 
this notation and terminology is common in commutative group theory,
non-commutative
group theorists often call $\bigoplus_{i\in I} G_i$ the {\em direct product\/}
of the family 
$\{G_i:i\in I\}$ and use product notation $\prod_{i\in I} G_i$ instead of 
$\bigoplus_{i\in I} G_i$. Since 
the former 
notation could easily lead to confusion with {\em Cartesian
products\/}, especially among topologists, we decided to use
 the ``commutative looking'' notation $\bigoplus_{i\in I} G_i$ instead 
 of $\prod_{i\in I} G_i$ common in non-commutative group theory.
However, for a {\em finite\/} family of groups $G_1,G_2,\dots, G_n$ we will use
the product notation $G_1\times G_2\times \dots \times G_n$ instead of 
 $G_1\oplus G_2\oplus \dots \oplus G_n$.

Our next lemma exhibits a particular situation when the assumption of Corollary \ref{approximation:theorem} holds:
\begin{lemma}
\label{countable:summands}
Let $G=N\times \left(\bigoplus_{i\in I} G_i\right)$, where $N$ is an Abelian group and each group $G_i$ is countable. 
Then the family 
$$
\mathcal{C}=\left\{N'\times\left(\bigoplus_{j\in J} G_j\right): J\in[I]^{\le\omega} \mbox{ and }N'\in\mathcal{S}(N)\right\}
$$
is a club in $[G]^{\le\omega}$ such that $\mathcal{C}\subseteq \mathcal{N}_G$.
\end{lemma}

\begin{proof}
$\mathcal{C}$ is trivially a club in $[G]^{\le\omega}$.  Let us see that $N'\times\left(\bigoplus_{j\in J} G_j\right)$ is Markov 
embedded into $G$ for every $J\in [I]^{\le\omega}$ and each countable subgroup $N'$ of $N$. As a direct 
summand, the subgroup $\bigoplus_{j\in J} G_j$ of $\bigoplus_{i\in I} G_i$ is super-normal in $\bigoplus_{i\in I} G_i$.
The subgroup $N'$ of the Abelian group $N$ is trivially super-normal  in $N$. 
This implies that $N'\times\left(\bigoplus_{j\in J} G_j\right)$ is super-normal in $N\times \left(\bigoplus_{i\in I} G_i\right)=G$.
From Corollary \ref{super-notmal:implies:topemb} and 
Lemma \ref{topological:embedding:implies:markov:embedding}, we now conclude that 
$N'\times\left(\bigoplus_{j\in J} G_j\right)$ is Markov embedded in $G$.
\end{proof}

\begin{corollary}
\label{direct:sum:corollary}
If $G=N\times\left(\bigoplus_{i\in I} G_i\right)$, where $N$ is an Abelian group and each group $G_i$ is countable, then Markov 
and Zariski topologies on $G$ coincide.
\end{corollary}

\begin{proof} Apply
Corollary \ref{approximation:theorem} and Lemma \ref{countable:summands}.
\end{proof}

\begin{corollary}
\label{direct:sum:corollary2}
If $G=\bigoplus_{i\in I} G_i$, where each group $G_i$ is countable,  then Markov and Zariski topologies on $G$ coincide.
\end{corollary}

In our  last lemma we offer a formal extension of the last corollary to certain subgroups of direct sums of countable groups.

\begin{lemma}\label{co:summands}
Let $H$ be a subgroup of $G= \bigoplus_{i\in I} G_i$, where each group $G_i$ is countable. For every $J\in [I]^{\le\omega}$ define $G_J= \bigoplus_{i\in J} G_i$
and $H_J=H\cap G_J$. If the family  $\mathcal{J}=\{J\in [I]^{\le\omega}: H_J\in  \mathcal{N}_{G_J}\}$ contains
a club in $[I]^{\le\omega}$, then Markov and Zariski topologies on $H$ coincide.
\end{lemma} 

\begin{proof}
Fix $J\in\mathcal{J}$. 
Note that $G=G_J\times G_{I\setminus J}$ and $H_J=G_J\cap H$ is Markov 
embedded in $G_J$ by the definition of $\mathcal{J}$. Applying 
Lemma \ref{(*)} we conclude that $H_J$ is Markov embedded in $H$ as well.
Thus $H_J\in\mathcal{N}_H$.

Let $\mathcal{C}\subseteq \mathcal{J}$ be a club in $[I]^{\le\omega}$. Consider  the map $\theta:[I]^{\le\omega}\to [H]^{\le\omega}$ defined by $\theta(J)=H_J$ for every $J\in[I]^{\le\omega}$. Note that $\theta$ is monotone, i.e., 
$J,J'\in [I]^{\le\omega}$ and $J\subseteq J'$ implies $H_J=\theta(J)\subseteq
\theta(J')=H_{J'}$. From this one can easily conclude 
that $\{\theta(J):J\in\mathcal{C}\}=\{H_J:J\in\mathcal{C}\}$ is a club in 
$[H]^{\le\omega}$. The conclusion of our lemma now follows from Corollary \ref{approximation:theorem}.
\end{proof}

\section{Connections with non-topologizable groups}

\begin{definition}\label{Def:topologizable} 
{\rm  Recall that a group $G$ is said to be {\em non-topologizable} if the only 
Hausdorff group topology of $G$ is the discrete one. 
A group $G$ is {\em topologizable\/} if it admits a non-discrete Hausdorff group topology.
}
\end{definition}

\begin{lemma}
\label{non-topologizable:lemma}
\begin{itemize}
\item[(i)] $G$ is non-topologizable if and only if  $\Mar_G$ is discrete.
\item[(ii)] If $\Zar_G$ is discrete, then $G$ is non-topologizable.
\end{itemize}
\end{lemma}
\begin{proof}
Item (i) is obvious. Item (ii) follows from Fact \ref{Markov:stronger:than:Zariski} and item (i). 
\end{proof}

The following lemma is easy to check.
\begin{lemma}\label{discrete:Zariski} 
The Zariski topology $\Zar_G$ of a group $G$ is discrete if and only  if there exist elementary algebraic sets $E_1, \ldots ,E_n$ such that 
$E_1\cup \ldots \cup E_n=G\setminus \{e_G\}$. 
\end{lemma}

The problem to construct a (countable) non-topologizable group was  raised by Markov and resolved consistently in \cite{Shelah} (see more 
details in Remark \ref{Shelah's_example} (iii)).  In  \cite{O} Ol$'$shanskij  used Lemma \ref{non-topologizable:lemma}(ii) and \ref{discrete:Zariski} 
 to produce the first ZFC solution of Markov's problem on the existence of non-topologizable 
countable groups (Ol$'$shanskij used an appropriate quotient of the (countable) Adian group $A(n,m)$). 

Now we give a sufficient condition (due to Shelah) that ensures that an uncountable group is  non-topologizable.

\begin{proposition}\label{discrete:Markov} 
{\rm (Shelah)}
An uncountable group $G$ is non-topologizable whenever the following two conditions hold: 
\begin{itemize}
\item[(a)] there exists $m\in \N$ such that $A^m=G$ for every subset $A$ of  $G$ with $|A|=|G|$; 
\item[(b)] for every subgroup $H$ of $G$ with $|H|<|G|$ there exist  $n\in \N$ and  
$x_1, \ldots, x_n\in G$ such that the intersection $\bigcap_{i=1}^n x_i^{-1}Hx_i$ is  finite. 
\end{itemize}
\end{proposition}

\begin{proof} Let $\mathcal{T}$ be a Hausdorff group topology on $G$.  There exists a $\mathcal{T}$-neighbourhood $V$  of $e_G$ with $V\ne G$.
Choose a $\mathcal{T}$-neighbourhood $W$ of $e_G$ with $W^m\subseteq V$.
 Now $V\ne G$ and  (a) yield $|W|<|G|$.  Let $H=\langle W \rangle$. Then $|H|=|W|\cdot\omega<|G|$. By (b) 
the intersection $O=\bigcap_{i=1}^n x_i^{-1}Hx_i$ is 
finite for some $n\in \N$ and elements $x_1, \ldots, x_n\in G$. Since each $x_i^{-1}Hx_i$ 
is a $\mathcal{T}$-neighbourhood of $e_G$, this proves  
that $e_G\in O\in\mathcal{T}$. Since $\mathcal{T}$ is Hausdorff, it follows that $\{e_G\}$ is $\mathcal{T}$-open,
and therefore $\mathcal{T}$ is discrete.
\end{proof}

\begin{remark}\label{Shelah's_example}
{\rm
\begin{itemize}
\item[(i)]
Note that in item (b) the number $n$ may depend of $H$, while in  item (a) the number $m$ is {\em the same} for all $A\in [G]^{|G|}$. 
(Indeed, one can easily see that in the circle group $G=\T$, written  additively, every neighbourhood 
$A$ of 0 in the usual topology satisfies $mA=G$ for some $m$ {\em depending 
on} $A$.) 
\item[(ii)]
Even the weaker form of (a) (with $m$ depending on $A\in  [G]^{|G|}$), yields that every proper subgroup of $G$ has size $<|G|$
(in the case $|G|=\omega_1$, the groups with this property are known 
as {\em Kurosh groups}, the first consistent example of a Kurosh group was given in \cite{Shelah}). 
\item[(iii)]
The above criterion was used by {Shelah}  \cite{Shelah} to produce the  first consistent example of a non-topologizable group (he worked under 
the assumption of CH and produced a group $G$ of size $\omega_1$ satisfying (a) with $m=10000$ and (b) with $n=2$). 
\end{itemize}}
\end{remark}
 
 \begin{lemma}\label{Examples1} 
 \begin{itemize}
  \item[(a)] If $G$ is non-topologizable, then the Markov embedded subgroups 
of $G$ are non-topologizable as well. In particular, no infinite Abelian subgroup of $G$ is  Markov embedded in $G$. More specifically, 
in a torsion-free non-topologizable group $G$ all cyclic subgroups are  not Markov embedded. 
\item[(b)] Every non-topologizable group $H$ is \topemb\ (and thus Markov embedded) in any ambient group $G$. 
  \end{itemize}
\end{lemma}

\begin{proof} (a) $\Mar_G$ is discrete by Lemma \ref{non-topologizable:lemma}(i).
If $H$ is a Markov embedded subgroup of $G$, then  
$\Mar_H=\Mar_G\restr_H$ must be discrete as well.
Applying Lemma \ref{non-topologizable:lemma}(i) once again, we conclude that
$H$ is non-topologizable.

(b) Let $H$ be a non-topologizable subgroup of a group $G$.
If $\mathcal{T}$ is a Hausdorff group topology on $H$, then $\mathcal{T}$ 
must be discrete, and so we can trivially extend $\mathcal{T}$  by
taking the discrete topology on $G$.  
\end{proof}

\begin{lemma}
\label{absolute:markov}
For a group $H$ the following conditions are equivalent:
 \begin{itemize}
   \item[(i)]  $H$ is Markov embedded in every group $G$ that contains it as a subgroup,
   \item[(ii)] $H$ is non-topologizable.
\end{itemize}
\end{lemma}
\begin{proof}
(i)$\to$(ii). According to \cite{T}, $H$ admits an embedding into some non-topologizable group $G$. Since $\Mar_G$ is discrete, from (i) we conclude that $\Mar_H$ must also be discrete. Hence $H$ is non-topologizable.

(ii)$\to$(i) follows from Lemma \ref{Examples1}(b) and Theorem \ref{topological:embedding:implies:markov:embedding}.
\end{proof}

\begin{remark}\label{absolutely-Markov:implies:absolutely-Hausdorff?}
{\rm 
\begin{itemize}
\item[(i)]
Lemma  \ref{absolute:markov} should be compared to Theorem \ref{Main:theorem}.
According to Lemma  \ref{topological:embedding:implies:markov:embedding} 
 for countable groups $H$ Theorem \ref{Main:theorem}  along with 
 Corollary \ref{Mar:Embedded>Zariski:Embedded:corollary}
provides many examples of countable Abelian groups that are Markov embedded in every subgroup 
where they are embedded as a {\em normal} subgroup. 
\item[(ii)]
According to Lemma \ref{absolute:markov} no infinite Abelian group $H$ can be 
Markov embedded in every group containing $H$ as a subgroup (as infinite Abelian groups are topologizable).
\end{itemize}
}
\end{remark}

According to \cite{KT} there exists a countable torsion-free non-topologizable group $G$. Hence for this group no cyclic
 subgroup $C$ is Markov embedded into $G$ by Lemma \ref{Examples1}. 
The next proposition yields that none of them is a {\em normal\/} subgroup of $G$.

\begin{proposition}\label{cyclic:subgroups}
If a group $G$ has an infinite cyclic subgroup as a normal subgroup, 
then $G$ is topologizable. 
\end{proposition}

\begin{proof}
Assume $H$ is an infinite normal cyclic subgroup of $G$. Then $H$ is Hausdorff embedded in $G$
by Corollary \ref{cyclic:subgroup}. So $H$ is Markov embedded in $G$ as well (Lemma \ref{topological:embedding:implies:markov:embedding}). Since $\Mar_H$ is non-discrete, it follows that $\Mar_G$ is non-discrete as well. Hence $G$ is topologizable. 
\end{proof}

\begin{lemma}
\label{topologization:criterion:using:Zariski:topology}
Let $G$ be a countable group for which the Zariski topology $\Zar_G$ is not discrete. Then $G$ is topologizable.
\end{lemma}

\begin{proof}
Since $\Zar_G=\Mar_G$ by Fact \ref{Markov:theorem}, from our assumption
it follows that $\Mar_G$ is not discrete. Now apply Lemma \ref{non-topologizable:lemma}(i). 
\end{proof}

\begin{corollary}
\label{topologizable:corollary}
Let $G$ be an infinite countable group such that $\Zar_G$ is compact. Then $G$ is topologizable.
\end{corollary}

\begin{proof}
An infinite compact space cannot be discrete, and the result follows from 
Lemma
\ref{topologization:criterion:using:Zariski:topology}.
\end{proof}

\begin{remark}
\label{(a*)fails}
{\rm 
\begin{itemize}
   \item[(i)] According to Remark \ref{Remark:ZvsM}, item (b) of Lemma \ref{Markov:vs:Zariski:embedding}
has the following stronger (but non-symmetric) form: (b$^*$) {\em  If $\Zar_G=\Mar_G$ and $H$ is  Zariski embedded in $G$, then $H$ is also Markov embedded in $G$   and $\Zar_H=\Mar_H$.} 

\item[(ii)]  According to \cite{S1} there exists a countable non-topologizable subgroup $H$ of a non-topologizable group $G$
with $\Mar_G \ne \Zar_G$. Then $\Mar_H$ is the discrete topology of $H$, 
thus $H$ is Markov embedded in $G$ by item (b) of Example \ref{Examples1}.
Since $H$ is countable, $\Zar_H=\Mar_H$ (Fact \ref{Markov:theorem}), and so 
the topology $\Zar_H$ is discrete. Therefore, $H$ is also Zariski embedded in $G$ by 
item (a) of Lemma \ref {Markov:embedding:coincides:with:Zariski:embedding:for;countable:groups}.
We see that $H$ is both Markov and Zariski embedded in $G$, and yet $\Zar_G\ne \Mar_G$. Therefore, the stronger
form (a$^*$) of item (a) of Lemma \ref{Markov:vs:Zariski:embedding} obtained by adding the condition $\Zar_G = \Mar_G$ to 
the conclusion of (a) may fail.

\item[(iii)] Item (i) and Remark \ref{Remark:ZvsM} explain why we preferred to announce Lemma 
\ref{Markov:vs:Zariski:embedding} in its present form that gives a pleasing symmetry between  items (a) and (b) of this lemma. 
\end{itemize}
}
\end{remark}

\section{Absolutely \topembing s and a gap in \cite{S2}}
\label{error:section}

Recall that an Abelian group $G$ is called {\em indecomposable\/} if for every direct 
product 
decomposition $G=G'\times G''$ either $G'=\{0\}$ or $G''=\{0\}$.

\begin{definition} 
{\rm
\cite{DS1}
A group $G$ is called:
\begin{itemize}
\item[(a)]  {\em absolutely \topemb\/} provided that $G$ is \topemb\ in every group $H$ containing $G$ as a normal subgroup.
\item[(b)]  {\em absolutely Zariski embedded\/} provided that $G$ is Zariski embedded in every group $H$ containing $G$ as a normal subgroup.
\end{itemize}
}
\end{definition}

By Lemma \ref{topological:embedding:implies:markov:embedding} and  Corollary \ref{Mar:Embedded>Zariski:Embedded:corollary} 
every absolutely \topemb\ Abelian group is also absolutely Zariski embedded. Even if the latter property may seem weaker,
one can show that it imposes a very strong restraint on the structure of an Abelian group. 

\begin{theorem}\label{Main:theorem} \cite{DS1}
Every absolutely Zariski embedded (in particular, every absolutely \topemb)
Abelian group is indecomposable.
\end{theorem} 

\begin{remark}\label{Example:nonZar/Mar_embedded}
{\rm 
Let $H$ be a decomposable Abelian group. By
Theorem \ref{Main:theorem},
there exists a group $G$ containing $H$ as a (normal) subgroup 
such that 
$H$ is not Zariski embedded in $G$.  
By Corollary \ref{Mar:Embedded>Zariski:Embedded:corollary}, $H$ 
is not Markov embedded in $G$ either.
}
\end{remark} 

The next characterization obtained in \cite{DS1} shows, among other things,
that divisible Abelian groups are never absolutely \topemb:
 
\begin{theorem}\label{2main:theorem}\cite{DS1}
An Abelian group $G$ is  absolutely \topemb\ if and only if 
the identity map $id_G$ of $G$ and minus the identity map $-id_G$ of $G$ are the only automorphisms of $G$. 
\end{theorem} 

In \cite{S2} the author takes a subgroup $H$ of a direct product $G=\prod_{\alpha\in I} G_\alpha$ of countable groups 
$G_\alpha$, then finds a countable 
set $I^*\subseteq I$ and considers the countable normal subgroup $G^*=H\cap \prod_{\alpha\in I} G^*_\alpha$
of $H$, where $G^*_\alpha=G_\alpha$ for $\alpha\in I^*$ and $G^*_\alpha=\{1_\alpha\}$
for $\alpha\in I\setminus I^*$. (Here $1_\alpha$ denotes the identity element of 
$G_\alpha$.) Then a metric group topology $\mathcal{T}^*$ on $G^*$ is constructed, and the author says: ``{\em Since $G^*$ is a normal subgroup of $H$, the neighbourhoods of the identity in the topology $\mathcal{T}^*$ form a neighbourhood base
of the identity of some group topology $\mathcal{T}$ on $H$.\/}'' The author than employs this topology $\mathcal{T}$ on $H$ to finish the proof of the main result of \cite{S2}: $\Mar_H=\Zar_H$ for such an $H$. (Apparently, the product $\prod_{\alpha\in I} G_\alpha$ in the author's terminology is what we call a direct sum
$\bigoplus_{\alpha\in I} G_\alpha$.)

The italicized statement above is an {\em essential gap\/} in the proof of the main result of \cite{S2}, and so the result 
itself should be considered an unsolved open problem (see our Question \ref{subgroups:of:products:question}).
Indeed, for this proof to work one has to extend a Hausdorff group topology $\mathcal{T}^*$ on $G^*$ to some Hausdorff 
group topology $\mathcal{T}$ on a bigger group $H$ containing $G^*$ as a normal subgroup.
Since {\em a priori\/} there is no control whatsoever either over the topology $\mathcal{T}^*$ on 
$G^*$ or the ambient group $H$ this 
appears impossible unless the group $G^*$ in question is \topemb\ in {\em every\/} group $H$ containing $G^*$ as a normal subgroup.
That is, $G^*$ apparently must be absolutely \topemb.
Indeed, our next lemma clearly demonstrates the inherent non-triviality of this extension problem.

\begin{lemma}
\label{metric:witness}
Let $G^*$ be a countable group that is not absolutely \topemb. Then there exists a group $H^*$ containing $G^*$ as a normal subgroup of countable index  and a metric group topology $\mathcal{T}^*$
on $G^*$ that cannot be extended to any Hausdorff group topology on $H^*$.
\end{lemma}

\begin{proof}
Indeed, by our assumption and a result from \cite{DS1} there exists a group $H^*$ containing $G^*$ as a normal subgroup of 
countable index such that $G^*$ is not \topemb\ in $H^*$.
By  Theorem  \ref{metric:groups}, there must exist a metric group topology $\mathcal{T}^*$ on $G^*$ that cannot be 
extended to any Hausdorff group topology on $H^*$. \end{proof}

From this lemma we ought to conclude that the proof of the main result of \cite{S2} could possibly work 
only in the case when $G^*$ is absolutely \topemb. The class of (countable) absolutely
 \topemb\ groups is extremely narrow. Indeed, according to Theorems \ref{2main:theorem} and \ref{Main:theorem}, an 
Abelian absolutely \topemb\ group $G^*$ must be indecomposable and every automorphism of $G^*$ must be either the identity map of $G^*$
or minus the identity map of $G^*$.

A persistent reader might still feel that there are yet additional circumstances in the setting of \cite{S2} that are not 
accounted for in Lemma \ref{metric:witness}. Indeed, the group $G^*$ in question appears to be ``nicely embedded'' in the direct product. 
However, Lemma \ref{lemma:of:product:of:two} demonstrates that this is a mere illusion. From this lemma one has to conclude the following: Given 
a subgroup $H$ of the square $G'\times G'$ of a countable group $G'$, one cannot reasonably expect 
to be able to extend a metric group topology from the (normal) subgroup $G^*=H\cap(G'\times\{1_{G'}\})$ of $H$ to any Hausdorff group topology on $H$ 
unless $G^*$ does not admit a decomposition $G^*=G_1\times G_2$ such that both groups $G_1$ and $G_2$ are infinite. Moreover, if 
$G^*$ has such a decomposition,
then $G^*$  may even fail to be both Markov embedded in $H$ and Zariski embedded in $H$.

\section{Final remarks and open questions}

Let $\mathcal{MZ}$ be the class of groups $G$ for which Markov and Zariski  topologies coincide: $\Mar_G=\Zar_G$. 

It might be tempting to generalize Corollary \ref{direct:sum:corollary2}
even further:

\begin{question}
\label{subgroups:of:products:question}
 Let $G$ be a subgroup of a direct sum of countable groups. Does $G$ belong to $\mathcal{MZ}$?
\end{question}

An attempt of providing a positive answer to  this question has been recently made in \cite{S2} but the proof contains essential errors that were pointed out in 
Section \ref{error:section}.

Since every Abelian group is a subgroup of a direct sum of countable groups,
a positive answer to Question \ref{subgroups:of:products:question} would yield that
Corollaries \ref{direct:sum:corollary} and \ref{direct:sum:corollary2}
are equivalent. A rather limited partial positive answer to Question
 \ref{subgroups:of:products:question} can be found in Lemma \ref{co:summands}.

\begin{question}\label{Ques4}
  \begin{itemize}
    \item[(i)]  Is  $\mathcal{MZ}$ closed under finite direct sums?
    \item[(ii)] Is $\mathcal{MZ}$ closed under arbitrary direct sums?
  \end{itemize}
\end{question}

\begin{question}\label{Ques5}
Let $H$ be an Abelian subgroup of $G$.
\begin{itemize}
  \item[(i)] If index of $H$ in $G$ is finite, does $G$ belong to $\mathcal{MZ}$?
  \item[(ii)] If index of $H$ in $G$ is countable, does $G$ belong to $\mathcal{MZ}$?
\end{itemize}
What is the answer to both items (i) and (ii) if one additionally assumes that $H$ is a normal subgroup of $G$?
\end{question}

As witnessed by Lemma \ref{Markov:vs:Zariski:embedding}, 
the following question is ultimately related to the Markov's problem:

\begin{question}
Let $H$ be a (normal) subgroup of a group $G$.
\begin{itemize}
  \item[(i)] If $H$ is Markov embedded in $G$, must $H$ also be Zariski embedded in $G$?
  \item[(ii)] If $H$ is Zariski embedded in $G$, must it also be Markov embedded in $G$?
\end{itemize}
\end{question}

Our next two questions should be compared with Lemma \ref{absolute:markov}.
\begin{question}
Describe the class of groups $G$ such that $G$ is Markov embedded in every group that contains $G$ as a normal subgroup.
\end{question}

Note that an Abelian group $G$ with the above property must be absolutely Zariski embedded by 
Corollary \ref{Mar:Embedded>Zariski:Embedded:corollary}, hence $G$ must be indecomposable by Theorem \ref{Main:theorem}.

\begin{question}
Let $H$ be a group that is Zariski embedded in every group $G$ containing $H$ as a subgroup. Must $H$ be non-topologizable?
\end{question}

Let $\{G_i:i\in I\}$ be a family of groups, and for every $i\in I$, let 
$H_i$ be a \topemb\ subgroup of $G_i$. Then the Cartesian product $\prod_{i\in I} H_i$ 
is a \topemb\ subgroup of the Cartesian product $\prod_{i\in I} G_i$. 
Similarly, the direct sum $\bigoplus_{i\in I} H_i$ is a \topemb\ subgroup 
of the direct sum $\bigoplus_{i\in I} G_i$.
In other words, the class $\mathcal{H}$ of \topembing s is closed under both Cartesian products and direct sums.
We do not know whether the remaining two classes $\mathcal{M}$ of Markov embeddings and $\ZZ$ of Zariski embeddings 
are closed under taking Cartesian 
products and direct sums. In fact, even in the weakest possible form, this is an open question: 

\begin{question}\label{Question:productivity}
Assume $H_1$ is Markov (Zariski) embedded in $G_1$. Is it true that 
for every group $G_2$ the subgroup $H_1\times G_2$ of $G=G_1\times G_2$ 
is  Markov (Zariski) embedded in $G$? 
\end{question}

\begin{remark}{\rm If 
Question \ref{Question:productivity} has a positive answer 
in the particular case when $H_1$ is a Zariski embedded subgroup
of an Abelian group $G_1$, then the more general Corollary \ref{direct:sum:corollary} would follow from 
the less general Corollary \ref{direct:sum:corollary2}.
Indeed, let $D$ be any divisible Abelian group containing $N$ (for example, the divisible 
hull of $N$). Then $D$ is a direct sum of countable groups 
\cite{Fuchs}. 
Therefore, for the group $\widetilde{G}= D\times\left(\bigoplus_{i\in I} G_i\right)$ 
one has $\Mar_{\widetilde{G}}=\Zar_{\widetilde{G}}$ by Corollary \ref{direct:sum:corollary2}. 
Since $D$ is Abelian, $N$ is super-normal in $D$ (Corollary \ref{Cor(b)}), and thus $N$ is Markov embedded in $D$ by Corollary \ref{super-notmal:implies:topemb} and Lemma \ref{topological:embedding:implies:markov:embedding}.
By Corollary \ref{Mar:Embedded>Zariski:Embedded:corollary}, $N$ is also Zariski embedded in $D$.
According to the presumed positive answer of Question \ref{Question:productivity} in this case,
$G$ now becomes Zariski embedded in $\widetilde{G}$. Since $\Mar_{\widetilde{G}}=\Zar_{\widetilde{G}}$, Remark 
\ref{Remark:ZvsM} would finally yield $\Mar_G=\Zar_G$. }
\end{remark}

\end{document}